\documentclass[11pt,a4paper,reqno]{amsproc}
\usepackage{eucal}

\usepackage{amsmath}
\usepackage{amsthm}
\usepackage{amssymb}

\newtheorem{thm}{Theorem}[section]

\newtheorem{lem}[thm]{Lemma}

\theoremstyle{remark}

\numberwithin{equation}{section}

\makeatletter
\newcommand{\puteqnnum}{\refstepcounter{equation}\hfill\@eqnnum}
\newcommand{\eqnnumlbl}[1]{\puteqnnum\label{#1}}

\def\resetschnoerkel{\def\schnoerkel{\relax}}%
\resetschnoerkel
\newcounter{saveeqn}

\let\theequationoriginal=\theequation

\def\thearabicequation{\ifnum\value{section}=0 \arabic{equation}\schnoerkel%
    \else\arabic{section}.\arabic{equation}\schnoerkel\fi}
\def\thealphequation{\ifnum\value{section}=0
            \mbox{\arabic{saveeqn}\alph{equation}\schnoerkel}
       \else \mbox{\arabic{section}.\arabic{saveeqn}\alph{equation}\schnoerkel}\fi}
\def\theequation{\thearabicequation}

\newcommand{\alpheqn}[1][\relax]{
     \refstepcounter{equation}
     \if#1\relax \relax
       \else \label{#1}
     \fi  
     \setcounter{saveeqn}{\value{equation}}%
    \setcounter{equation}{0}%
    \renewcommand{\theequation}{\thealphequation}}
\newcommand{\reseteqn}{\setcounter{equation}{\value{saveeqn}}%
     \renewcommand{\theequation}{\thearabicequation}}

\newcommand{\myref}[1]{{\normalfont (\ref{#1})}}

\newcommand{\plref}[1]{{\normalfont \ref{#1}}}

\newlength{\blockwidth} 
\newlength{\blockmargin} 
\newlength{\blockrightmargin} 
\newcommand{\numblock}[1]{\par\medbreak\noindent\hspace*{\blockmargin}
   \parbox{\blockwidth}{#1}%
    \stepcounter{equation}\hfill\@eqnnum}

\newcommand\ga{\alpha} 
  
\newcommand\go{\omega}

\newcommand{\C}{\mathbb{C}}
\newcommand{\R}{\mathbb{R}}

\newcommand{\Z}{\mathbb{Z}}

\newcommand\ca{\mathcal{A}}

\newcommand\cd{\mathcal{D}}

\newcommand\ch{\mathcal{H}}

\newcommand\cl{\mathcal{L}}

\newcommand\im{\operatorname{im}}

\newcommand\ind{\operatorname{ind}}

\newcommand\loc{\operatorname{loc}}

\newcommand\Op{\operatorname{Op}}

\newcommand\sign{\operatorname{sign}}
\newcommand\spec{\operatorname{spec}}

\newcommand\tr{\operatorname{tr}}

\newcommand\vol{\operatorname{vol}}

\newcommand{\cinfz}[1]{C_0^\infty(#1)}

\newcommand\LIM{\operatorname*{LIM}}

\newcommand\DST{\displaystyle}

\newcommand{\comment}[1]{\relax}

\newcommand{\ovl}{\overline}

\newcommand{\restr}{\restriction}

\newcommand{\class}{\operatorname{cl}}

\newcommand{\ASS}{AvrSeiSim:IPP}
\renewcommand{\tilde}{\widetilde}

\newcommand{\OpAn}{{{\rm Op}^0(A)}}

\renewcommand{\Op}{{\rm Op}}

\newcommand{\specess}{\operatorname{spec}_{\operatorname{ess}}}

\newcommand{\scalar}[2]{\langle #1,#2\rangle}

\theoremstyle{plain}

\newcommand{\bigsetdef}[2]{\bigl\{ #1 \,\bigm|\, #2\bigr\}}

\newlength{\boxwidth}
\newenvironment{MLnumlist}{
\setlength{\boxwidth}{\textwidth}
\addtolength{\boxwidth}{-\leftmargin}
\addtolength{\boxwidth}{-1.5cm}
\newcommand{\MLitem}[1]{\item  \parbox[t]{\boxwidth}{##1}}}{}

\newcommand{\specflow}{\operatorname{sf}}

\makeatother

\begin{document}
\setcounter{page}{1}

\title[Spectral theory of boundary value problems]{Spectral 
theory of boundary value problems for Dirac type operators}

\author{Jochen Br\"uning}
\address{Humboldt Universit\"at zu Berlin\\
Institut f\"ur Mathematik\\ 
Unter den Linden 6\\
D--10099 Berlin
}
\curraddr{}
\email{bruening@mathematik.hu-berlin.de}

\author{Matthias Lesch}
\address{Humboldt Universit\"at zu Berlin\\
Institut f\"ur Mathematik\\ 
Unter den Linden 6\\
D--10099 Berlin
}
\curraddr{Universit\"at Bonn\\ Mathematisches Institut\\ Beringstr. 1\\
D--53115 Bonn}
\email{lesch@mathematik.hu-berlin.de, lesch@math.uni-bonn.de}
\thanks{Both authors were supported in part by Deutsche Forschungsgemeinschaft}

\date{}
\thanks{This paper is in final form and no version of it will
be submitted for publication elsewhere.}

\begin{abstract}
The purpose of this note is to describe a unified approach to the 
fundamental results in the spectral theory of boundary value problems, 
restricted to the case of Dirac type operators. 
Even though many facts are known and well presented in the literature 
(cf. the monograph of Booss-Wojciechowski \cite{BooWoj:EBPDO}), 
we simplify and extend or sharpen most results by using systematically 
the simple structure which Dirac type operators display near the boundary. 
Thus our approach is basically functional analytic, and consequently we 
achieve results which apply to more general situations than compact 
manifolds with boundary.

The details of the material presented here will be published elsewhere.

Thanks are due to R. T. Seeley and the referee for helpful comments.
\end{abstract}

\maketitle

\section{The compact case}
Let $\tilde{M}$ be a compact oriented Riemannian manifold of dimension $m$, 
$\tilde{E} \to \tilde{M}$ a hermitian vector bundle, 
and $\tilde{D}$ a symmetric elliptic differential operator of order 
$d\in\Z_+$ on $C^\infty (\tilde{E})$. The classical results on existence, 
uniqueness,  and regularity of solutions of $\tilde{D}$ are among the 
cornerstones of Global Analysis:
\begin{itemize}
\begin{MLnumlist}
\MLitem{$ \tilde{D}$ is essentially self-adjoint in $ L^2(\tilde{E})$ 
with domain $C^\infty(\tilde{E})$ 
(by slight abuse of notation, we denote the closure of 
$\tilde{D}$ by the same symbol);}\eqnnumlbl{intro-1.1}
\MLitem{
$\tilde{D}$ is a Fredholm operator (of index $0$),
i.e. there are a bounded operator $\tilde{Q}$ and compact operators 
$\tilde{K}_r , \tilde{K}_l$ in $L^2(\tilde{E})$ such that}\eqnnumlbl{intro-1.2}
\[
 \tilde{D}\tilde{Q} = I - \tilde{K}_r,\qquad \tilde{Q}\tilde{D} = I-\tilde{K}_l ;
\]
\MLitem{with respect to the Sobolev scale 
$H_s(\tilde{E}):= {\cd}((\tilde{D}^2 + I)^{s/2d}),$ $s\ge 0,$ $\tilde{Q}$
is of order $-d$  and $\tilde{K}_{r/l} $ of order $- \infty.$}
\eqnnumlbl{intro-1.3}
\end{MLnumlist}
\end{itemize}

The restriction to symmetric operators is not essential since we may 
always consider a given elliptic operator together with its adjoint. 
But it is a technical advantage for more refined questions like 
index theorems: We bring in an isometric involution, $\tilde\omega,$ 
on $\tilde E$ which anticommutes with $\tilde{D}$ on $C^\infty (\tilde{E})$ 
and hence produces a splitting
\begin{equation}
 \tilde{D} = \begin{pmatrix} 0 & \tilde D_- \\ \tilde D_+ & 0\end{pmatrix}
 \text{ on } \; C^\infty (\tilde{E}_+) \oplus  C^\infty (\tilde{E}_-),
  \label{intro-1.4}
\end{equation}
with $\tilde{E}_\pm$ the $\pm 1$--eigenbundle of $\tilde{\omega}$. 
Then, by the well--known formula of McKean-Singer we have
\begin{equation}
 \ind \tilde{D}_+ = \tr \bigl[ \tilde{\omega} e^{-t\tilde{D}^2}\bigr],
  \quad t>0.\label{intro-1.5}
\end{equation}
It is hence of great importance that 
(cf. eg. \cite[Lemma 1.9.1]{Gil:ITHEASIT2E})
\begin{itemize}
\begin{MLnumlist}
\MLitem{
for any differential operator, $\tilde{P}$, of order $p$ on 
$C^\infty(\tilde{E})
$, we have an asymptotic expansion}
\begin{equation}
\tr \Bigl[ \tilde{P} e^{-t\tilde{D}^2}\Bigr] \sim_{t\to 0+} 
  \sum a_j (\tilde{D},\tilde{P}) t^{(j-m-p)/d}.
  \label{intro-1.6}
\end{equation}
\end{MLnumlist}
\end{itemize}

Even though it took a long time after the original proof of the 
Atiyah--Singer Index Theorem \cite{AtiSin:IEOCM}, \cite{Pal:SASIT}
until a complete proof could be based on \myref{intro-1.4} 
and \myref{intro-1.5} (cf. \cite{AtiBotPat:HEIT}, \cite{Gil:CELEC}), 
the heat equation method now seems to be the most powerful tool for 
extensions of the Index Theorem. Recall that for the important class 
of twisted Dirac operators, with $\tilde{E} = \tilde{S}\oplus\tilde{F}, 
\tilde{S}$ a spin bundle on $\tilde{M}$, this theorem reads
\begin{itemize}
\item $\DST\ind \tilde{D}_+ = \int\limits_{\tilde{M}} \hat{A} (\tilde{M}) 
\wedge \rm ch \,\tilde{F}.$\eqnnumlbl{intro-1.7}
\end{itemize}

\section{Compact manifolds with boundary}
We want to present the extension of the main 
results quoted above (\myref{intro-1.1}, \myref{intro-1.2}, 
\myref{intro-1.3}, \myref{intro-1.6}, \myref{intro-1.7}) 
to Dirac type operators on manifolds with boundary. Though a good part 
of our results is more or less known, we obtain a conceptually as well as 
technically transparent derivation of this theory, with considerable 
simplifications and extensions in most cases. Moreover, the functional 
analytic approach we have developed lends 
itself naturally to substantial generalizations, e.g. to situations 
with non-compact boundaries.

To explain our work in greater detail, we consider a compact hypersurface,
$N$, in $\tilde{M}$  which bounds an open subset, $M$, of $\tilde{M}$. 
We assume that $N$ is oriented as the boundary of $M$.
Then we put 
$E:=\tilde{E} \restr M,$ \linebreak[3]
$D:=\tilde{D} \restr C^\infty (\overline{M})$. \linebreak[3]
Then $D$ is a first order elliptic differential operator on $M$, 
and symmetric in $L^2(E)$ with domain $C^\infty_0(E)$. If $D$ is of 
\emph{Dirac type} (in the sense of
\cite[Chap. 4, 4.4]{Gil:ITHEASIT2E})
then we obtain in a tubular neighbourhood, $U$, of $N$ 
in $\tilde{M}$ a very simple separation of variables. In fact, $U$ is 
isometric to $(-\varepsilon_0, \varepsilon_0) \times N$ with metric 
$dx^2\oplus g_N(x), x\in (-\varepsilon_0, \varepsilon_0),$
$g_N$ a smooth family of metrics on $N$, and with 
$E_N := \tilde{E} \restr N$ we obtain the following result.
\begin{lem}\label{intro-S1} 
Let $D$ be of Dirac type. As operator in $L^2(E\restr U)$ with domain 
$C^\infty_0(E\restr U)$, $D$ is unitarily equivalent to an operator of the form
\begin{equation}
\gamma\left(\frac{d}{dx} +A(x)\right)
\label{intro-1.8}
\end{equation}
in $L^2((-\varepsilon_0,\varepsilon_0),L^2(E_N))$
with domain $C^\infty_0((-\varepsilon_0,\varepsilon_0),
C^\infty(E_N)),$
where $\gamma\in{\cl}(L^2(E_N))$ and $A(x)$ is (the closure of) 
a symmetric elliptic differential operator on $E_N$ of first order; 
${\cd}(A(x))=: {\cd}$ is independent of $x$ and $A(x)$ depends smoothly 
on $x\in (-\varepsilon_0,\varepsilon_0).$

Moreover, the following relations hold:
\begin{subequations}\label{intro-1.9}
\begin{align}
  &\gamma^\ast =-\gamma,\quad \gamma^2=-I,\label{intro-1.9a}\\
  &\gamma({\cd})= {\cd} \text{ and } \gamma A(x) + A(x)\gamma = 0, 
\quad x\in(-\varepsilon_0,\varepsilon_0).\label{intro-1.9b}
\end{align}
\end{subequations}
\end{lem}
This lemma has been widely used for some time, especially in the product case 
$(g_N(x)\equiv g_N(0))$ where it plays a prominent role in \cite{AtiPatSin:SARG}. 
For non-product metrics some care is needed to compute $A(x)$ in each specific
case,
cf. e.g. \cite[Sec. 3.10]{Gil:ITHEASIT2E} and 
\cite[Sec. 5]{Bru:LITCMROT}.

We will base our analysis on a thorough study of the operator equation 
\myref{intro-1.8} with the structure properties \myref{intro-1.9};
\emph{these properties will be assumed throughout this paper}.
This is reasonable since the results we are aiming at can be obtained from 
merging "interior analysis" (to be carried out on $\tilde{M}$) with 
"boundary" analysis involving the operator \myref{intro-1.8}.

The main difference between the analysis of $D$ and $\tilde{D}$ lies, of
course, in the fact that $D$ is not  essentially self-adjoint on
$C^\infty_0(E).$ 
Moreover, if self-adjoint extensions of $D$ exist, they may differ widely with 
respect to existence, uniqueness, regularity, and heat trace expansions. 
It is, therefore, our first task to characterize those self-adjoint extensions 
which behave nicely with respect to existence, uniqueness, and regularity.

\section{Results for the model operator}
Replacing in \myref{intro-1.8} $L^2(E_N)$ by an arbitrary Hilbert space, $H$, 
and $C^\infty(E_N)$ by the domain, $H_1$, of  a self-adjoint operator $A$ in
$H$, we obtain the model operator
\begin{equation}
 D=\gamma\left(\frac{d}{dx} +A\right) \text{ in } \,{\ch}_0:=L^2(\R_+, H)
   \text{ with domain } C^\infty_0(\R_+, H_1).\label{intro-1.10}
\end{equation}

We will have to deal with variable coefficients but for the purpose of this
overview we will restrict to the constant coefficient case.
Indeed, for most of the problems dealt with in this paper operators with
variable coefficients merely appear as perturbations of constant ones,
in view of the Kato-Rellich Theorem.

On $C^\infty_0(\R_+, H_1)$ we clearly have
\begin{equation}
   (Df, g)_{{\ch}_0} - (f, Dg)_{{\ch}_0} = \scalar{f(0)}{\gamma g(0)}_H.
   \label{intro-1.11}
\end{equation}
Now if $D$ is symmetric on a subspace, ${\cd}^0$, of $C^\infty_0(\R_+, H_1)$
then it follows from \myref{intro-1.11} that, with $I-P$ the
orthogonal projection onto $\ovl{\cd^0}$ in $H$, we have
\begin{subequations}\label{intro-1.12}
\begin{equation}
 {\cd}^0 \subset {\cd}_P := \bigsetdef{f\in C^\infty_0 (\R_+, H_1)}{P f(0)=0}
  \label{intro-1.12a}
\end{equation}
and
\begin{equation}
  I - P \le \gamma^\ast P\gamma.\label{intro-1.12b}
\end{equation}\end{subequations}
Moreover, 
$D_P := D \restr {\cd}_P$ is a symmetric extension of 
$D\restr \cd^0$. If we assume for a moment that $H$ is of
finite dimension then it is readily seen that $D_P$ is essentially self-adjoint
in ${\ch}_0$ if and only if
\begin{equation}
  I - P = \gamma^\ast P\gamma.\label{intro-1.13}
\end{equation}
Indeed, if $D_{\max}$ denotes $(D\restr C^\infty_0 ((0,\infty), H_1))^\ast$ then
\begin{equation}
  {\cd}(D_{\max}) \subset H_{1,\loc} (\R_+, H),\label{intro-1.14}
\end{equation}
and \myref{intro-1.11} remains valid with $D_{\max}$ in place of $D$, 
for $f, g \in {\cd}(D_{\max}).$

An orthogonal projection, $P$, with \myref{intro-1.13} will be called 
$\gamma$--\emph{symmetric}; it is easy to see that $\gamma$--symmetric 
projections -- and hence self--adjoint extensions of $D$ -- exist 
if and only if 
\begin{equation}
  \sign (i \gamma \restr \ker A) = 0.\label{intro-1.15}
\end{equation}
As an illustration, note that $i\frac{d}{dx}$ does not admit self-adjoint 
extensions in $L^2(\R_+).$

Returning to the general case, we meet the essential difficulty that 
\myref{intro-1.14} has no reasonable analogue. In particular, elements of 
${\cd}(D_{\max})$ do not  admit $H$-valued restrictions to zero. To overcome
this obstacle, we imitate the Sobolev scales $H_s(E_N)$ and $H_s(E)$ and their
interplay in our abstract setting (which has some tradition in
Analysis, cf. eg. \cite[Chap. XIII]{Pal:SASIT}). $H_s(E_N)$ is replaced by
\begin{equation}\begin{split}
 H_s&:= H_s(A)\\
    & := \left\{\begin{array}{l}
  {\cd}(|A|^s), \text{ equipped with the graph norms for } s\ge 0; \\
\text{a suitable dual of } H_{-s} (A), \text{ for } s<0.\end{array}\right.
		\end{split}
\end{equation}
We also need
\begin{align*}
  H_\infty&:= H_\infty(A):=\bigcap_{s\in\R} H_s(A);\\
\intertext{and}
  H_{-\infty}&:=H_{-\infty}(A):=\bigcup_{s\in\R} H_s(A).
\end{align*}
Next we introduce, for $n\in\Z_+$,
\begin{subequations}\label{intro-1.17}
\begin{equation}
 {\ch}_n := {\ch}_n (\R_+, A) := \bigcap^n_{k=0} H_k (\R_+ , H_{n-k}(A)),
 \label{intro-1.17a}
\end{equation}
where, for $i,j\in\Z_+ ,$
\begin{equation}\begin{split}
  H_i(&\R_+, H_j (A)) \\
     &:= \bigsetdef{f\in{\ch}_0}%
{\Bigl(\frac{d}{dx}\Bigr)^l f\in L^2(\R_+, H_j(A)), \quad 0\le l\le i}.
  \label{intro-1.17b}
		\end{split}
\end{equation}
\end{subequations}
By duality and interpolation, we then obtain a scale of Hilbert spaces, 
${\ch}_s ={\ch}_s (\R_+, A), s\in \R.$ 
Generalizing the classical Trace Theorem for Sobolev spaces, 
we have the following result about trace maps which will allow 
the formulation of boundary conditions.
\begin{thm}\label{intro-S1a}
The map 
$$
 r : C^\infty_0 (\R_+, H_\infty)\ni f \longmapsto f(0) \in H_\infty 
$$
extends by continuity to a map
$$
  r_s : {\ch}_s \longrightarrow H_{s-1/2}, \quad s> 1/2,
$$
and also to a map
$$ 
 r^\ast : {\cd} (D_{\max}) \to H_{-1/2} .
$$
\end{thm}
Of course, the loss of regularity under the trace map requires the 
(continuous) extension of the boundary projections to the space $H_{-1/2}$. 
To deal with this, we introduce \emph{operators of finite order} on the 
Hilbert scale $(H_s(A))_{s\in\R}$. 
Thus, a linear map, $B:H_\infty\to H_\infty,$ is an operator of order
$\mu\in\R$ if for each $s\in\R$ there is a constant $C(s)$ such that, for any
$x\in H_\infty$,
\begin{equation}
 \| Bx\|_{H_s} \le C(s) \| x\|_{H_{s+\mu}} \,.\label{intro-1.18}
\end{equation}
In particular, $B$ extends to an element of ${\cl}(H_s, H_{s-\mu})$ 
for all $s\in\R.$ The totality of such operators forms the linear space 
$\Op^\mu(A).$ $\Op^{-\infty}(A) := \bigcap_{\mu\in\R}\Op^\mu (A)$ 
is called the space of \emph{smoothing operators}.

Thus we will have to require that the boundary projections are elements of 
$\Op^0(A)$. It follows easily from \myref{intro-1.18} that $\Op^0(A)$ is a 
$\ast$-algebra but it is, in general, not \emph{spectrally invariant} 
in the sense that $B\in \Op^0(A)$ and $B$ invertible in ${\cl}(H)$ 
implies $B^{-1}\in \Op^0(A)$.
To allow for a minimum of functional constructions, we do need actually 
even more. We are forced to restrict attention to certain subalgebras, 
$\Psi^0(A) \subset \Op^0(A),$ satisfying the following two conditions. 

\begin{itemize}
\item[($\Psi$1)] $\Psi^0(A)$ is a $\ast$-subalgebra of
$\OpAn$ with \emph{holomorphic functional calculus};\\
\item[($\Psi$2)] $\Psi^0(A)$ contains an orthogonal projection,
$P_+(A)$, satisfying 
\begin{equation}
 I-P_+(A)=\gamma^*P_+(A)
\gamma,\quad P_{(0,\infty)} (A)\le P_+(A)\le P_{[0,\infty)}(A).
 \label{intro-1.19}
\end{equation}
\end{itemize}

Recall that an algebra, ${\ca}\subset {\cl}(H),$ has holomorphic 
functional calculus if for $B\in {\ca}$ and $f$ holomorphic in a 
neighbourhood of $\spec B$ (in $\cl(H)$) we have $f(B) \in{\ca}$, 
where $f(B)$ is defined by the Cauchy integral. 
Note also that the existence of $P_+(A)$ with \myref{intro-1.19} 
is equivalent to 
\myref{intro-1.15} in the finite dimensional case.
In general, if $0 \notin \specess A$ then \myref{intro-1.15} is equivalent to 
the existence of a spectral projection of $A$ satisfying \myref{intro-1.19}.

Let us illustrate these conditions for the case where $A$ is an elliptic
differential 
operator on $C^\infty(E_N)$ and $N=\partial M$ as in Sec. 2.
Then we have $H_s(A) \simeq H_s(E_N)$, and a 
natural choice of the algebra $\Psi^0(A)$ is the algebra of classical 
pseudodifferential operators on $E_N$, to be denoted by 
$\Psi^0_{\class} (E_N)$. It follows from results of Seeley 
\cite[Thm. 5]{See:CPEO} that $\Psi^0_{\class}(E_N)$ has 
holomorphic functional 
calculus. Moreover, since $0 \notin \specess A$, we have to verify 
\myref{intro-1.15} to obtain a spectral projection, $P_+(A)$, of $A$ 
fulfilling \myref{intro-1.19}; but this is a consequence of the 
Cobordism Theorem.
To see this, 
we split $H=: H_+\oplus H_-$ according to the $\pm i$--eigenspaces of
$\gamma$. In view of 
\myref{intro-1.9},
\begin{equation}
  A=\begin{pmatrix} 0 & A_- \\ A_+ & 0 \end{pmatrix},
   \label{I.G1-1.17}
\end{equation}
and $A_+$ is a Fredholm operator with index
\begin{equation}
 \ind A_+ = \dim \ker A\cap \ker(\gamma-i) -\dim \ker A \cap \ker (\gamma + i).
 \label{I.G1-1.18}
\end{equation}
Now it is straightforward to check that there exists 
an orthogonal projection $P_+(A) \in \OpAn$ with the property
\myref{intro-1.19} if and only if
\begin{equation}
  \ind A_+ =0,
  \label{I.G1-1.19}
\end{equation}
and this follows from the Cobordism Theorem
(cf. the discussion after \cite[Cor. 3.6]{BruLes:EICNLBVP}). 

Again from Seeley's work, we deduce that 
$P_+(A) \in\Psi_{\class}^0(E_N)$ so ($\Psi$1), ($\Psi$2) are satisfied 
in a natural way.

Now we are in the position to formulate our results for the model operator. 
The main theorem concerning regularity reads as follows.
\begin{thm}\label{intro-S2}
Let $H$ be a Hilbert space and $A$ a self-adjoint operator in $H$. 
Assume, moreover, that an algebra, $\Psi^0(A)\subset \Op^0(A),$ 
is given with the properties \textnormal{($\Psi$1)} and \textnormal{($\Psi$2)}.

Then $D_P$ with domain \myref{intro-1.12a} is essentially self-adjoint in 
$L^2(\R_+, H)$ for any orthogonal projection $P\in \Psi^0(A)$ with the 
properties
\begin{equation}
 \gamma^\ast P\gamma = I - P\tag{\ref{intro-1.13}}
\end{equation}
and
\begin{equation}
  (P,P_+(A)) \text{ is a Fredholm pair.}\label{intro-1.20}
\end{equation}
The domain of the closure of $D_P$ is 
$\cd(D_P)=\bigsetdef{f\in\ch_1(\R_+,A)}{Pf(0)=0}$.

Conversely, if $A$ is discrete then the self-adjointness of $D_P$
on $\bigsetdef{f\in\ch_1(\R_+,A)}{Pf(0)=0}$
implies \myref{intro-1.13} and \myref{intro-1.20}.
\end{thm}

We note that the orthogonal projection $P_-(A)=I-P_+(A)\in\Psi^0(A)$
obviously does not satisfy \myref{intro-1.20}. However, it can be
shown that $D_{P_-(A)}$ is essentially self--adjoint with domain
$\cd(D_{P_-(A)})\supsetneqq\bigsetdef{f\in\ch_1(\R_+,A)}{Pf(0)=0}$.
Hence the "self--adjointness" in the last statement of the Theorem
cannot be replaced by "essentially self-adjoint on $\cd_P$".

Recall that a pair of (orthogonal) projections, $(P_1,P_2)$, in $H$
is said to form a \emph{Fredholm pair} if the map
\[P_2:P_1(H)\longrightarrow P_2(H)\]
is Fredholm. This is easily seen to be equivalent to the fact that
\begin{subequations}\label{intro2-1.21}
\begin{align}
   &P_2P_1(H) \text{ (and hence } P_1P_2(H)) \text{ is closed in } H\\
\intertext{and}
   & \dim(\ker P_2\cap \im P_1)+\dim (\ker P_1\cap \im P_2)<\infty
\end{align}
\end{subequations}
(such that $(P_1,P_2)$ is Fredholm if and only if so is $(P_2,P_1)$).
In this case, one calls
\begin{equation}
    \ind(P_1,P_2):=\dim(\ker P_2\cap \im P_1)-\dim(\ker P_1\cap \im P_2)
  \label{intro-1.22}
\end{equation}
the \emph{index of the pair} $(P_1,P_2)$; cf. \cite[IV.4.1]{Kat:PTLO},
\cite{Boj:ALCPFPS}, \cite[Sec. 24]{BooWoj:EBPDO}, \cite{\ASS}.

The proof of Theorem \plref{intro-S2} is technically complicated and cannot
be described here in detail. We have to interpolate 
various abstract concepts of regularity such that
we can show their mutual equivalence step by step.

We can view Theorem \plref{intro-S2} as the analogue of \myref{intro-1.1}
for the model operator. Taking advantage of the self-adjointness of $D_P$ 
we can try to satisfy \myref{intro-1.2} by setting
\begin{equation}
  Q := \int\limits_{|\lambda|\ge 1} \lambda^{-1} d  E(\lambda),
  \label{intro-1.21}
\end{equation}
where $E(\lambda) = E_{D_P}(\lambda), \lambda\in \R,$ denotes the spectral 
resolution of $D_P$. 
From our abstract regularity with respect to the Sobolev scale
$\ch_s(\R_+,A)$, $s\in\R$, and a standard compactness argument we then derive
the following analogue
of \myref{intro-1.2} and \myref{intro-1.3}.
\begin{thm}\label{intro-S3}
We assume the situation of Theorem \plref{intro-S2} and, in addition, 
that $A$ is discrete. For $\phi\in C^\infty_0(\R)$ with $\phi=1$ near 
$0$ we put $Q_\phi := \phi  Q$. Then $Q_\phi$ maps into ${\cd}(D_P)$ 
and there are compact operators, $K_{r/l,\phi}$, 
in ${\ch}$ such that
\begin {equation}
  D_PQ_\phi = \phi - K_{r,\phi},\qquad 
  Q_\phi D_P = \phi -K_{l,\phi}.
  \label{intro2-1.22}
\end{equation}
Moreover, $Q_\phi $ is of order $-1$ and $K_{r/l,\phi}$ of order $-\infty$ 
with respect to the Sobolev scale ${\ch}_s (\R_+, A), s\in\R.$
\end{thm}
We have pointed out that, for the purposes of regularity theory, we may treat
variable coefficients as a perturbation of the constant coefficient case;
the same is true for index theory, by deformation. As in 
\cite{BruLes:EICNLBVP}, a basic tool of our analysis is the analogue
of a formula due to Sommerfeld, expressing the heat kernel of the
model operator $D_{P_+(A)}=: D_{P_+}$ for $x,y,t>0$,
\begin{equation}\begin{split}
 e^{-t D_{P_+}^2}&(x,y)\\
  &=(4\pi t)^{-1/2} 
      \left(e^{-(x-y)^2/4t}+(I-2P_+)e^{-(x+y)^2/4t}\right) e^{-tA^2}
      \\
    &\quad+ (\pi t)^{-1/2}(I-P_+)
      \int_0^\infty e^{-(x+y+z)^2/4t} A e^{Az-t A^2} dz.
    \label{intro2-1.23}
		\end{split}
\end{equation}
In order to derive a reasonable index formula, we assume now that
$(A+i)^{-1}\in\cl_p(H)$ (the Schatten--von Neumann class) for some
$p>0$, and that
\begin{equation}
   \eta(A;s):=\frac{1}{\Gamma\bigl(\frac{s+1}{2}\bigr)}
  \int_0^\infty t^{(s-1)/2}\tr_H\bigl[Ae^{-tA^2}\bigr] dt
\end{equation}
extends to a meromorphic function in $\C$ without a pole at $0$.
Then it is not difficult to derive from \myref{intro2-1.23} the following
result which computes the "boundary contribution" to the index. Recall
that $\go:=\tilde\go\restr E$ is the involution defining the
index. Then $\go$ has to commute with $A$ and $P$, and we put
\begin{equation}
   A^+:=A \frac 12(\go+I),\quad P^+:=P\frac 12(\go+I).
\end{equation}
Furthermore, $\LIM\limits_{t\to 0+} f(t)$ denotes the constant term in the
asymptotic expansion of $f$, and $\eta(A^+):=\eta(A^+;0)$.

\begin{thm}\label{intro-S5} Under the assumptions of Theorem \plref{intro-S3},
let $P$ also commute with $\go$ and $A$,
\begin{equation}
   [P,A]=[P,\go]=0.
    \label{intro2-1.25}
\end{equation}
Then we have, for any $\phi\in\cinfz{\R}$ with $\phi=1$ in a neighborhood
of $0$,
\begin{equation}\begin{split}
   \LIM_{t\to 0+} \tr\bigl[\phi\go e^{-t D_P^2}\bigr]&=-\frac
   12\bigl(\eta(A^+)+\dim\ker A^+\bigr)+\ind(P_{\ge 0}(A^+),P^+)\\
    &=:\xi(A^+)+\ind(P_{\ge 0}(A^+),P^+).
		\end{split}
\end{equation}
\end{thm}
It turns out that regular boundary projections other than $P_+(A)$ are quite
useful. Moreover, their functional analytic characterization shows that they
form a convenient tool for "glueing indices". We have used this already
in \cite{BruLes:EICNLBVP} for the more complicated case of $\eta$--invariants;
we recall briefly what is involved.

Consider the model operator \myref{intro-1.10} on $L^2(\R,H)$
where it is essentially self--adjoint. The reflection isometry
(generated by $\sigma(x)=-x),$
\[\Phi_\sigma:L^2(\R,H)\longrightarrow L^2(\R_+,H\oplus H),\quad
   f\longmapsto (f\restr \R_+, f\circ \sigma\restr \R_+),
\]
transforms this operator unitarily to a model operator on $\R_+$,
\begin{subequations}
\begin{equation}
    D=\tilde \gamma\left(\frac{d}{dx}+\tilde A\right),
\end{equation}
with domain $\cinfz{\R_+,H_1\oplus H_1}$ in $L^2(\R_+,H\oplus H)$,
where
\begin{equation}
 \tilde \gamma=\begin{pmatrix}\gamma & 0\\ 0 & -\gamma\end{pmatrix}
 \quad\text{and}\quad \tilde A=\begin{pmatrix} A & 0 \\ 0 & -A\end{pmatrix}.
\end{equation}
The boundary condition (the smooth transmission condition) simply becomes,
for $u=(u_1,u_2)\in\cinfz{\R_+,H_1\oplus H_1},$
\begin{equation}
    u_1(0)=u_2(0).\label{intro2-1.27c}
\end{equation}
\end{subequations}
The subspace of $H\oplus H$ defined by \myref{intro2-1.27c} can
be viewed as the $+1$--eigenspace of the involution
\begin{equation}
   \tau:=\begin{pmatrix}  0&I\\I&0\end{pmatrix},\label{intro2-1.28}
\end{equation}
and if $P$ is any orthogonal projection in $H\oplus H$ satisfying
\begin{equation}
     \tau P=(I-P)\tau,
\end{equation}
then
\begin{subequations}\label{intro2-1.30}
\begin{equation}
      \ga:= P- (I-P)
\end{equation}
is a second self--adjoint involution with
\begin{equation}
    \tau\ga=-\ga \tau.
\end{equation}
\end{subequations}
Now we define a family of unitaries by
\begin{subequations}\label{intro2-1.31}
\begin{equation}
    U(\theta):=\cos \theta I + \sin \theta \ga\tau,\quad |\theta|<\pi/2,
\end{equation}
and a corresponding family of projections by
\begin{equation}\begin{split}
    P(\theta)&:= U(\theta)^*PU(\theta)\\
      &=\cos^2\theta P+\sin^2\theta (I-P) +\sin \theta \cos\theta \tau.
		\end{split}\label{intro2-1.31b}
\end{equation}\end{subequations}
Clearly, $P(\pi/4)$ projects onto the $+1$--eigenspace of $\tau$, so the
family \myref{intro2-1.31b} is a deformation of \myref{intro2-1.27c}.
In particular, from $\tau\tilde A+\tilde A\tau=0$ we see that
\begin{equation}
   P=P_+(\tilde A):=\begin{pmatrix}  P_+(A)&0\\0& I-P_+(A)\end{pmatrix}
\end{equation}
is an admissible choice, and we obtain a continuous deformation 
from the Atiyah--Patodi--Singer boundary condition - 
which decouples $D$ - to the smooth transmission boundary condition. 
Therefore, the following result is of interest for index calculations.

\begin{thm}\label{intro-S6}
If $P$ is regular, so is $P(\theta)$ for $|\theta|<\pi/2,$ and the family 
$D_{P(\theta)}$ is graph continuous in $(-\pi/2, \pi/2)$. 
\end{thm}
Index theorems for the model operator on finite or infinite intervals 
in $\R$ have also been used frequently in recent years, 
cf. eg. \cite{RobSal:SFMI}. The methods described so far are the 
starting point for a systematic abstract study of this topic. 
Thus, let us consider the operator
\begin{equation}
  D^{a,b}_{P,Q} := \frac{d}{dx} + A(x)
\end{equation}
with domain (cf. \myref{intro-1.17a})
${\cd}^{a,b}_{P,Q} := \bigsetdef{f\in \ch_1([a,b], A)}%
{Pf(a) = (I-Q)f(b)=0}$ in $L^2([a,b], H)$. 
Now it is of importance that we admit variable coefficients. 
Precisely, we assume that $A\in C^\infty ([a,b], {\cl}(H_1,H))$ 
and that each $A(x)$ is self-adjoint in $H$ with domain $H_1$, 
and a Fredholm operator. 
It is well known that the index of $D^{a,b}_{P,Q}$ is related to the 
spectral flow of the family $A(x)$ across $[a,b],$ to be denoted by 
$\specflow_{[a,b]} A$ (cf. part II of \cite{AtiPatSin:SARG} for the definition 
and the main properties of $\specflow A$). 
Our main result here has various obvious 
corollaries like a formula for the spectral flow and for the variation 
of $\xi$-invariants; it reads as follows.
\begin{thm}\label{intro-S7}
Let $A$ be as above and assume that $P,Q$ are orthogonal projections in 
$H$ such that $(P_{\ge 0}(A(a)),P)$ and $(P_{\ge 0}(A(b)),Q)$ are 
Fredholm pairs. 
Then (the closure of) $D^{a,b}_{P,Q}$ is a Fredholm operator with
\begin{equation}\begin{split}
 \ind D^{a,b}_{P,Q} &= -\frac{1}{\sqrt{\pi}} \int\limits^b_a 
   \LIM\limits_{t\to 0+} t^{1/2} \tr_H \bigl[A'(x) e^{-tA(x)^2}\bigr] dx\\
  &\quad -\xi(A(a)) + \ind(P_{\ge 0} (A(a)),P)\\
  &\quad + \xi(A(b)) - \ind (P_{\ge 0} (A(b)), Q)\\
 &= \specflow_{[a,b]} A + \ind (P_{\ge 0} (A(a)),P) - 
     \ind (P_{\ge 0} (A(b)), Q).
		\end{split}
\end{equation}
\end{thm}

Finally, we turn to the asymptotic expansion of the heat trace. 
Again, we obtain by simple abstract methods a substantial result, 
without using any pseudodifferential technology but, for the time being, 
our result is far from being best possible. Thus, we now restrict 
attention to projections, $P$, which satisfy the assumptions in Theorem 
\plref{intro-S2} and, in addition,
\begin{equation}
  P-P_+ (A) \in \Op^{-1} (A).
  \label{intro2-1.35}
\end{equation}
Now, for the expansion problem variable coefficients are considerably 
more difficult than constant ones. We thus allow a family, 
$(A(x))_{x\ge 0}$, satisfying the conditions needed in 
Theorem \plref{intro-S7} and also, with $A:=A(0)$,
\begin{equation}
  A(x) \equiv A \quad \text{ for }\; x\ge x_0 >0.
  \label{intro-1.36}
\end{equation}
Since we need the expansion result only near the boundary 
$x=0$ in most applications, \myref{intro-1.36} does not mean an 
essential restriction, but \myref{intro2-1.35} does.

To obtain complete asymptotic expansions, we need of course to 
assume expansion properties for the family $(A(x))$. 
These assumptions are somewhat technical so we refrain from 
stating them here in detail, it may suffice to assert that 
they are satisfied in the most prominent case when 
$\Psi^0 (A) = \Psi^0_{\class} (E_N),$ for some compact Riemannian 
manifold $N$ and some hermitian bundle $E_N\to N$. Then we obtain
\begin{thm}\label{intro-S8}
Under the assumptions mentioned above, we obtain for any 
$\phi\in C^\infty_0 (\R_+)$ an asymptotic expansion
\begin{equation}
 \tr_{L^2(\R_+,H)} \bigl[\phi e^{-tD^2_P}\bigr] 
 \sim_{t\to 0_+} \sum_{\substack{j\ge 0\\ 0\le k\le 1}} 
   a_{jk}(P,\phi)t^{j/2-p}\; \log^k t,
\end{equation}
for some $p>0$ depending on $A$. If $j/2 -p\le 0$ we have $a_{j1}=0.$
\end{thm}

\section{Results for manifolds with boundary}
It is not difficult to translate Theorems \plref{intro-S2} through
\plref{intro-S7}
into statements on $D$ and the Sobolev scales $H_s(E)$ and 
$H_s(E_N), s\in \R$. First, we only need to make 
Lemma \plref{intro-S1} somewhat 
more explicit. For this, we introduce, on $U$, the global coordinate
\begin{equation}
 x(p) := \operatorname{dist}(p,N) , \quad p\in N,
 \label{intro-1.23}
\end{equation}
and denote by
\begin{equation}
 \Phi : L^2(E\restr U) \longrightarrow 
  L^2 ((-\varepsilon_0, \varepsilon_0), L^2(E_N))
  \label{intro-1.24}
\end{equation}
the isometry implicit in Lemma \plref{intro-S1}. Then we have the properties
\begin{subequations}\label{intro-1.25}
\begin{align}
&\psi\Phi u= \Phi \bigl((\psi\circ x)u\bigr),
   &&\psi \in C^\infty (-\varepsilon_0, \varepsilon_0), u\in L^2(E\restr U),
              \label{intro-1.25a}\\ 
    & (\Phi u)(0) = u \restr N, & &u\in C^\infty (E\restr \overline{M}),
               \label{intro-1.25b}\\
   &  \Phi\bigl((\psi\circ x) H_s(E)\bigr) =  
      \psi{\ch}_s, 
   &&\psi\in C^\infty_0 (-\varepsilon_0,\varepsilon_0), s\in \R,
\end{align}
\end{subequations}
which allow us to localize near $N$ and to transfer regularity.

To formulate the boundary conditions, we restrict attention to orthogonal 
projections in $L^2(E_N)$ which are classical pseudodifferential operators 
i.e. from now on we choose 
$\Psi^0(A) =\Psi^0_{\class} (E_N),$ as indicated above. 
In the theory of boundary value problems for linear elliptic differential 
operators, it was observed by Calder\'on \cite{Cal:BVPEE} 
that a prominent role is played by an idempotent, 
$C^+ \in\Psi^0_{\class} (E_N)$, with the property that
\begin{equation}\begin{split}
 C&^+(H_s(E_N)) = N_s(E_N)\\
  & := \bigsetdef{u\in H_s(E_N)}%
{u=\tilde{u}\restr N \text{ for } \tilde{u}\in H_{s+1/2}(E)
\text{ with } D\tilde{u} = 0}\label{intro-1.26}
		\end{split}
\end{equation}
for all $s\in\R; C^+ $ is called the \emph{Calder\'{o}n projector} 
(cf. \cite{See:SIBVP,See:TPO}, a comprehensive summary can be found in
\cite[Appendix]{Gru:HTEGDOWBC}). 
One checks that
\begin{equation}
 C^+ - P_+(A) \in \Psi^{-1}_{\class}(E_N)
\end{equation}
which explains the importance of the Atiyah--Patodi--Singer boundary condition. 
In order to obtain boundary conditions which define Fredholm operators 
(as in \myref{intro-1.2}, \myref{intro-1.3}), Seeley introduced the notion 
of "well-posed" boundary condition \cite{See:TPO}. 
Combining the results described so far with a microlocal argument,  
we obtain the following optimal version of Seeley's result, as a consequence 
of our general theory.
\begin{thm}\label{intro-S4}
Let $P\in\Psi^0_{\class} (E_N)$ be an orthogonal projection in 
$L^2(E_N)$ satisfying \myref{intro-1.13}.

Then $D$ with domain 
${\cd}_P := \bigsetdef{f\in H_1(\tilde{E}\restr\overline{M})}%
{P(f\restr N) = 0}$ is self-adjoint in $L^2(E)$ if and only if 
$P$ is well-posed in the sense of Seeley. This, in turn, is equivalent
to the fact that $(P,P_+(A))$ is a Fredholm pair.

In this case, denote by $D_P$ the closure of $D\restr \cd_P$. Then 
there are a bounded operator, $Q$, and compact operators, 
$K_{r}, K_{l},$ in $L^2(E)$ such that 
$Q$ maps into ${\cd}(D_P)$ and 
$$
  D_PQ = I-K_r, \quad QD_P= I-K_l.
$$
Moreover, with respect to the Sobolev scale $H_s(E), s\in \R, Q$ 
is of order $-1$ and $K_{r/l}$ of order $-\infty$.
\end{thm}
For index theorems in the framework of 
Theorem \plref{intro-S8}, we can reduce the computations to the 
situation addressed in Theorem \plref{intro-S5}, 
in view of the deformation properties of Fredholm pairs in 
$\Psi^0_{\class} (E_N)$ (here ($\Psi$1) is important) and the 
Local Index Theorem. Recall that the latter asserts that
$$
 \lim_{t\to 0+} \tr_{\tilde{E}} 
\bigl[\tilde{\omega} e^{-t\tilde{D}^2} (p,p)\bigr] \vol_M (p) =: \alpha_D (p)
$$
exists for any $p\in\tilde{M}$ and coincides with the Atiyah--Singer integrand.
We obtain the following result which contains the 
Atiyah--Patodi--Singer Theorem as well as a result of Agranovich and Dynin 
\cite{AgrDyn:GBVPESDD}.
\begin{thm}
Under the assumptions of Theorem \plref{intro-S8} we obtain
\begin{equation}
 \ind D_{P,+} = \int\limits_M \alpha_D + \int\limits_N \beta_D - 
    \xi(A(0)^+) + \ind (P_{\ge 0} (A(0)^+), P^+).
    \label{intro-1.38}
\end{equation}
Here, $(A(x))_{x\ge 0}$ is the family obtained from $\tilde{D}$ 
by separation of variables as in Lemma \plref{intro-S1}, 
and $\beta_D$ is given by a universal expression in the 
derivatives $A^{(j)} (0)$ .
\end{thm}
Geometric formulas for $\beta_D$ (for specific $\tilde{D}$) 
have been given by Gilkey \cite[Chap. 3,10]{Gil:ITHEASIT2E}; they involve 
the second fundamental form of $N$.

Theorem \plref{intro-S6} has also a very useful application in the 
present context to a "glueing formula" for indices. 
This is meant to say that we emphasize the constancy of the 
index along the curve given by $D_{P(\theta),+}$, for 
$ 0\le \theta \le \pi/4,$ decoupling the problem at $\theta = 0$.

More precisely, in the situation described at the beginning of 
Sec. 2, $M$ is decomposed as the union of two manifolds with 
boundary which we denote by $M^+$ and $M^-$,
\begin{equation}
 M=M^+ \cup M^- \cup N .
\end{equation}
Here, the orientation of $N$ is as boundary of $M^+$, 
hence the opposite of the boundary orientation induced by $M^-.$ 
The isomophism \myref{intro-1.24} in Lemma \plref{intro-S1} 
is normalized in such a way that
\begin{equation}
  \Phi (L^2(\tilde{E}\restr U\cap M^+)) = L^2 ((0,\varepsilon_0), L^2(E_N)).
  \label{intro2-1.40}
\end{equation}
Then we have the following result relating the index of $\tilde{D}$ 
to the index of
boundary value problems on $M^+$ and $M^-$ 
(which we add in the notation of the operators, for clarity). 
\begin{thm}\label{intro-S9}
Let $P\in\Psi^0_{\class} (E_N)$ satisfy \myref{intro-1.13} 
and \myref{intro-1.20}, and assume that $P$ commutes with $\omega$. Then
\begin{equation}
 \ind \tilde{D}_+ = \ind D^{M^+}_{P,+} + \ind D^{M^-}_{I-P,+}.
 \label{intro-1.41}
\end{equation}
\end{thm}
Let us elaborate a little bit on the formula \myref{intro-1.41}. 
To $D^{M^+}_{P,+} $ we can directly apply the index formula 
\myref{intro-1.38} which makes (implicitly) use of the chosen 
orientation through Lemma \plref{intro-S1}, 
whereas for $D^{M^-}_{I-P,+}$ we have to change the orientation of $N$. 
In \myref{intro-1.38}, this amounts to changing the sign of 
$\beta_D$ and replacing $A(0)^+ $ by $-A(0)^+$. 
Now an easy computation gives
\begin{align*}
  \ind (P_{\ge 0}(-A(0)^+),I-P^+)&= \ind (I-P_{>0} (A(0)^+), I-P^+)\\
    &= -\ind (P_{>0}(A(0)^+), P^+)
\intertext{and}
 \ind(P_{\ge 0} (A(0)^+),P^+)& = \ind (P_{\ge 0}(A(0)^+), P_{>0} (A(0)^+))\\
     &\quad  + \ind (P_{>0} (A(0)^+), P^+)\\
 & = \dim \ker A(0)^+ + \ind (P_{>0}(A(0)^+), P^+) .
\end{align*}
Since $\eta(-A)=-\eta(A),$ we see that \myref{intro-1.41} and 
\myref{intro-1.38} combine to
$$
  \ind \tilde{D}_+ = \int\limits_M \alpha_D ,
$$
as it must be.

We remark that  \myref{intro-1.41} remains true also in more 
general situations where $\tilde{M}$ need not be compact 
(but, of course, $\tilde{D}$ has to be Fredholm).

Theorem \plref{intro-S8} and \plref{intro-S9} can be used to give very 
simple and transparent proofs of various known results, 
notably the Cobordism Theorem and index theorems of Callias type. 
As mentioned before, the model operator does not "see" compactness; 
as an example of its wider applicability, we mention a very simple proof 
of the covering space version of the Atiyah--Patodi--Singer Theorem, 
due to Ramachandran.

Finally, we turn to the asymptotic expansion of the heat trace, 
based on Theorem \plref{intro-S8}. Applied to the situation at hand, 
our result is not as strong as the recent expansion theorem 
proved by G. Grubb \cite{Gru:HTEGDOWBC} even though we employ a 
much more elementary technique. Nevertheless, let us state Grubb's 
result for completeness (our proof needs \myref{intro2-1.35} 
for the time being).\pagebreak[3]
\begin{thm}
Let $P\in \Psi^0_{\class} (E_N)$ be well-posed and satisfy 
\myref{intro-1.13}. Then there is an asymptotic expansion
\begin{equation}
  \tr \bigl[ \omega^i D^j_P e^{-t D^2_P}\bigr] 
 \sim_{t\to 0+}\sum_{\substack{k\ge 0\\ l=0,1}} a^{ij}_{kl} (D,P)
  t^{(k-m-j)/2} \log^l t .
\end{equation}
Here,
\begin{equation}
  a^{ij}_{k1}= 0 \quad \text{ if } \, k-m-j< 0.\label{intro-1.43}
\end{equation}
\end{thm}
Under the more restrictive assumption \myref{intro2-1.35} 
we can show that \myref{intro-1.43} holds even for 
$k-m-j\le 0;$ this is relevant for the study of the 
($\zeta$-regularized) determinant of $D^2_P$ and $D_P.$


\end{document}